\documentclass[pdflatex,sn-mathphys-num]{sn-jnl}
\usepackage{lipsum}
\usepackage{amsfonts, amssymb}
\usepackage{amsthm}
\usepackage{graphicx}
\usepackage{epstopdf}
\usepackage{algorithmic}
\usepackage{mathtools}
\usepackage{cleveref}
\usepackage{todonotes}
\usepackage{subcaption}
\usepackage{comment}
\usepackage{bbm}
\usepackage{tikz}
\usepackage{xcolor}
\usepackage{gensymb}
\usepackage{macros}
\hypersetup{hidelinks=true}

\begin{document}

\title[Extensions for Fractional Sobolev Spaces on Lipschitz Submanifolds]{Extension Operators for Fractional Sobolev Spaces on Lipschitz Submanifolds}

\author*[1]{\fnm{Philipp} \sur{Weder}}\email{philipp.weder@epfl.ch}
\affil*[1]{\orgdiv{Institute of Mathematics}, \orgname{Ecole Polytechnique Fédérale de Lausanne}, \orgaddress{\street{Station 8}, \city{Lausanne}, \postcode{1015}, \state{Vaud}, \country{Switzerland}}}


\abstract{
It is well-known that any Lipschitz domain is an extension domain for $W^{s, p}$.
This paper extends this result to Lipschitz subsets of compact Lipschitz submanifolds of $\R^n$.
We adapt the construction of an extension operator for Lipschitz domains in \cite{di_nezza_hitchhikers_2012} to manifolds via local coordinate charts.
Moreover, the dependence of the continuity constant on the size of the extension domain is made explicit. This explicit scaling is crucial for applications in numerical analysis, such as a posteriori error estimation for geometry simplifications.}

\maketitle

\section{Introduction}
Sobolev extension operators find applications in numerical analysis, for instance finite elements in time-dependent domains~\cite{lehrenfeld_eulerian_2019}, homogenization \linebreak methods~\cite{hopker_note_2014}, or \textit{a posteriori} defeaturing error estimation~\cite{buffa_analysis-aware_2022}. For numerical applications, not only the existence of such operators but also the dependence of their continuity constants on the geometry of the extension domain are of practical relevance. In particular, in geometry simplification problems, the dependence of these constants on the size of the domain is crucial.

For a domain $\domain \subset \R^n, n \geq 2$, it is classical that under suitable regularity conditions on the boundary $\partial \domain$, there exists a bounded linear extension operator for the Sobolev spaces $\Wsp{\domain}$. Lipschitz domains provide a well-known class of extension domains for all fractional orders $s > 0$ ~\cite{morrey_lebesgue_1961, stein_singular_1970, mclean_strongly_2000, di_nezza_hitchhikers_2012}. A detailed proof for this statement for $s \in (0, 1)$ is provided by \citet{di_nezza_hitchhikers_2012}[§ 5] or for $s > 0$ by \citet{mclean_strongly_2000}. A general characterization of extension domains for the fractional orders $s \in (0, 1)$ and $p \in (0, \infty)$ is provided by \citet{zhou_fractional_2014}.

However, many applications require the study of Sobolev spaces on submanifolds of $\R^n$, for instance surfaces and interfaces arising in computational mechanics, geometry processing, and boundary integral methods. While Sobolev spaces on Lipschitz manifolds are well understood~\cite{hebey_nonlinear_2000}, extension results for Sobolev spaces defined on subsets of manifolds are less developed.
In this paper, we generalize the result for Lipschitz domains to subsets of compact Lipschitz submanifolds of $\R^n$, valid for all $p \in [1, +\infty)$ and $s \in (0, 1)$. The proof will follow the one in \cite{di_nezza_hitchhikers_2012}, but care must be taken when handling the local parameterizations of the manifold.

The remainder of the paper is structured as follows: First, we introduce the essential definitions for Sobolev spaces on manifolds. Next, we will present some preliminary lemmas, followed by the proof of the main statement. We will also outline an approach to generalize the result for Sobolev exponents $s \geq 1$. Finally, we will conclude with a summary of our findings.

\section{Preliminaries}
We first collect standard definitions and notation for Lipschitz domains, Lipschitz submanifolds, and fractional Sobolev spaces.

\begin{definition}[Lipschitz domain]
\label{def:lipschitz domain}
    A domain $\domain$ with boundary $\partial \domain$ is said to be a \emph{Lipschitz domain} if for every $x \in \partial \domain$, there exists a neighborhood $U(x)$ of $x$ and a Lipschitz-continuous function $\phi: \R^{n-1} \to \R$ such that, after relabeling and reorienting the coordinate axes if necessary, we have
    \begin{align*}
        \domain \cap U(x) = \{y \in U(x) \mid y_n > \phi(y_1, \ldots, y_{n-1})\}.
    \end{align*}
\end{definition}

In this paper, we are interested in Sobolev spaces on compact Lipschitz submanifolds, which naturally appear as boundaries of Lipschitz domains.

\begin{definition}[Lipschitz submanifold]
\label{def:lipschitz manifold}
Let $\manifold \subset \R^n$ be a closed bounded subset. Then, $\manifold$ is a bounded $k$-dimensional submanifold of $\R^n$ if for every point $x \in \manifold$, there exists an open set $U_x \subset \R^k$, an open neighborhood $V_x \subset \R^n$ of $x$ and a bi-Lipschitz homeomorphism $t_x: U_x \to V_x \cap \manifold$.
\end{definition}

We equip a Lipschitz submanifold $\manifold$ with the canonical surface measure induced by a parameterization, which we denote by $\dd{s}$.
Local computations are performed via the area formula. That is, let $(U_x, t_x)$ be a local parameterization as described above. We set
\begin{align*}
J_{t_x}(\hat{x}) \coloneqq 
\begin{cases}
    \sqrt{\det(Dt_x^\top(\hat{x}) Dt_x(\hat{x}))}, & t_x \text{ differentiable at } \hat{x},\\
    0, \text{ otherwise },
\end{cases}
\end{align*}
which is well-defined in $L^\infty(U_x)$, since $t_x$, being a bi-Lipschitz continuous bijection, is almost everywhere differentiable with respect to the Lebesgue measure due to Rademacher's Theorem \cite{leoni_first_2017}[Theorem 9.14]. Then, for any Borel-measurable function, $h: V_x \cap \manifold \to [-\infty, \infty]$, define
\begin{align}
\label{eq:change of variable formula}
\int_{V_x \cap \manifold} h(x) \dd{s}(x)  \coloneqq \int_{U_x} h(t(\hat{x})) J_{t_x}(\hat{x}) \dd{\hat{x}}.
\end{align}
Given an atlas $\{(U_i, t_i)\}_{i \in \mathcal{I}}$ of $\manifold$, where $\mathcal{I}$ denotes a countable index set, and a subordinate partition of unity $\{\phi_i\}_{i \in \mathcal{I}}$, the integral for a Borel-measurable function $h: M \to [-\infty, \infty]$ is defined as
\begin{align*}
    \int_M h(x) \dd{s}(x) \coloneqq \sum_{i \in \mathcal{I}} \int_{V_i \cap \manifold} h(x) \phi_i(x) \dd{s}(x).
\end{align*}
Note that the surface measure is independent of the choice of parameterization \cite{hebey_nonlinear_2000,lee_introduction_2012}. For subsets $\genbd \subset \manifold$, we denote by $|\genbd|$ the surface measure of the set.

We now state a preliminary lemma about the local parameterizations of a bounded Lipschitz submanifold. It ensures that the constants arising from local charts can be chosen uniformly over the manifold.

\begin{lemma}
\label{lem:local parameterizations}
    Let $\manifold$ be a bounded $k$-dimensional Lipschitz submanifold of $\R^n$. Then, the following assertions hold:

    \begin{enumerate}
        \item There exists $\epsilon > 0$, such that for all $x \in \manifold$, we can find a local parameterization $t_x: U_x \to \ball{\epsilon}{x} \cap \manifold$ with $0 \in U_x$ and $t_x(0) = x$.
        \item For $x \in \manifold$, denote by $L_{t_x}$ and $L_{t_x^{-1}}$ the Lipschitz constants of $t_x$ and $t_x^{-1}$, respectively.
        Then, there exist constants $L, K, J, I > 0$, such that for all $x\in \manifold$,
        \begin{align*}
            L_{t_x} &\leq L, & ||J_{t_{x}}||_{L^\infty(U_x)} \leq J,
            \\
            L_{t_x^{-1}} &\leq K, & ||J_{t_{x}^{-1}}||_{L^\infty(\ball{\epsilon}{x} \cap \manifold)} \leq I.
        \end{align*}
        \item For every $x \in \manifold$, we have $\ball{\epsilon / L}{0} \subset U_x \subset \ball{\epsilon K}{0}$.
    \end{enumerate}    
\end{lemma}

\begin{proof}
    Since $\manifold$ is bounded and therefore compact, we find a finite collection of bi-Lipschitz local parameterizations $s_j: U_j \to V_j \cap \manifold, j=1, \ldots, N$, for open sets $U_j \subset \R^k$ and $V_j \subset \R^n$, such that the collection $\{V_j \cap \manifold\}_{j=1}^N$ is an open cover of $\manifold$.

    For assertion 1, let $\epsilon > 0$ be the Lebesgue number of the open cover $\{V_j \cap \manifold\}_{j = 1}^N$. Let $x \in \manifold$ and $i$ the index such that $x \in V_i$. Then, we set $U_x \coloneqq s_i^{-1}(\ball{\epsilon}{x} \cap \manifold)$ and $t_x \coloneqq (s_i)_{|U_x}$. Translating $U_x$ by $s_i^{-1}(x)$ yields the assertion.

    For assertion 2, we denote by $L_j$ the Lipschitz constants of $s_j$, and we set $L \coloneqq \max_j L_j$. Then, we note that $t_x$ in the first step of the proof inherits the Lipschitz constant of $s_i$, which is not impacted by the translation. Hence, it follows immediately that $L_{t_x} \leq L$. The existence of $J$ follows from the fact the operator norm of $Dt_x$ is uniformly bounded by $L$ as well. The proof for the existence of $K$ and $I$ is analogous.

    For the last assertion, we fix $x \in \manifold$, and we choose $y, z \in \ball{\epsilon}{x} \cap \manifold$. Then, we have $|t_{x}^{-1}(y) - t_{x}^{-1}(z)| \leq K |y - z| < K \epsilon$. Conversely, for any $\hat{y} \in \ball{\epsilon / L}{0}$, we have $|t_x(\hat{y}) - x| \leq L |\hat{y}| \leq \epsilon$, which proves the assertion.
\end{proof}

\subsection{Function spaces on Lipschitz manifolds}

On a $k$-dimensional Lipschitz submanifold $\manifold$, we define the classical $\Lp{\manifold}$ space, for $p \in [1, \infty]$, with the norm
\begin{align*}
    \Lpnorm{u}{\manifold}^p &\coloneqq \int_\manifold |u(x)|^p \dd{s}(x) & \text{if } p < \infty,
    \\
    ||u||_{L^\infty(\manifold)} &\coloneqq \inf\{C \geq 0: |u| \leq C \text{ a. e. in } \manifold\} & \text{ if } p = \infty.
\end{align*}

Finally, we define the classical Sobolev spaces $\Wsp{\manifold}$ in terms of Gagliardo-Slobodeckij seminorms \cite{grisvard_elliptic_2011}: Let $s \in (0, 1)$ and $p \in [1, +\infty)$. Then, we define
\begin{align*}
    \Wsp{\manifold} \coloneqq \{u \in \Lp{\manifold}: \Wspnorm{u}{\manifold} < \infty\},
\end{align*}
where
\begin{align}
\label{eq:norm definition}
    \Wspnorm{u}{\manifold}^p &\coloneqq \Lpnorm{u}{\manifold}^p + \Wspseminorm{u}{\manifold}^p,\\ \Wspseminorm{u}{\manifold}^p &\coloneqq \int_\manifold \int_\manifold \frac{|u(x) - u(y)|^p}{|x - y|^{k + p s}}\dd{s}(x) \dd{s}(y).
\end{align}
The definition of Sobolev spaces on subsets $\genbd \subset \manifold$ with $|\genbd| > 0$ is analogous.

\begin{remark}
    Different choices of atlases yield equivalent norms for the Sobolev spaces $\Wsp{\manifold}$ \cite{mclean_strongly_2000}. Since $\manifold$ is fixed and compact, the equivalence constants depend only on the Lipschitz characteristics of the charts and are invariant under rescaling of the ambient space. Hence, the scaling relations in \cref{lemma:scaling properties} are independent of the choice of atlas.
\end{remark}

\section{The main result}
\label{sec:main result}
Let us now state the main result:

\begin{theorem}
\label{thm:boundary extension operator}
Let $\manifold$ be a $k$-dimensional Lipschitz submanifold of $\R^n$ according to \cref{def:lipschitz manifold} and $\genbd \subset \manifold$ an open connected subset such that its boundary $\partial \genbd$ is a $(k-1)$-dimensional Lipschitz submanifold.

Then, for $p \in [1, +\infty)$ and $s \in (0,1)$, there exists a linear and  continuous extension operator $\extop{\genbd}{\manifold}: \Wsp{\genbd} \to \Wsp{\manifold}$
such that
\begin{align*}
    \Wspnorm{\extop{\genbd}{\manifold}(u)}{\manifold}^p \leq C^{(0)}\left(1 + |\genbd|^{-\frac{sp}{k}} + |\genbd|^{\frac{(1-s)p}{k}}\right) \Lpnorm{u}{\genbd}^p +  C^{(s)}\Wspseminorm{u}{\genbd}^p,
\end{align*}
for constants $ C^{(0)} > 0$ and $ C^{(s)} > 0$ that depend on $n,p,s, \manifold$, and the Lipschitz regularity of the boundary $\partial \genbd$, but are independent of the measure $|\genbd|$.
\end{theorem}

Before providing a proof of \cref{thm:boundary extension operator}, we prove several preliminary lemmas and analyze the dependence of all appearing constants on the size of $\genbd$.
Using a partition of unity, we decompose the function into locally supported terms. This reduces the problem to extending functions that either vanish near $\partial \genbd$ (\cref{lemma:zero extension on boundaries}) or are supported within a single boundary chart (\cref{lemma:extension at boundaries}).

\begin{lemma}
\label{lemma:zero extension on boundaries}
    Let $\manifold$ be a $k$-dimensional Lipschitz submanifold of $\R^n$ and $\genbd \subset \manifold$ such that $|\genbd| > 0$.
    Let $u \in \Wsp{\genbd}$ with $s \in (0,1)$ and $p \in [1, +\infty)$. If there exists a compact subset $F \subset \genbd$ such that $u = 0$ in $\genbd \setminus F$, then the extension $u^\star$ defined by
    \begin{align*}
        u^\star \coloneqq \begin{cases}
            u(x), & x \in \genbd,\\
            0, & x \in \manifold \setminus \genbd,
        \end{cases}
    \end{align*}
    is an element of $\Wsp{\manifold}$, such that
    \begin{align*}
        \Lpnorm{u^\star}{\manifold} \leq \Lpnorm{u}{\genbd}
    &&
    \text{ and }
    &&
        \Wspseminorm{u^\star}{\manifold}^p \leq C^\star \Lpnorm{u}{\genbd}^p + \Wspseminorm{u}{\genbd}^p,
    \end{align*}
    with
    \begin{align*}
    C^\star \coloneqq \frac{|\manifold \setminus \genbd|}{\dist(\manifold \setminus \genbd, F)^{k + sp}}.
    \end{align*}
\end{lemma}

\begin{proof}
    Clearly, $u^\star \in \Lp{\manifold}$ with $\Lpnorm{u^\star}{\manifold} \leq \Lpnorm{u}{\genbd}$. Furthermore, we have
    \begin{align*}
        \Wspseminorm{u^\star}{\manifold}^p = \Wspseminorm{u}{\genbd}^p + 2 \int_{\manifold \setminus \genbd} \int_{\genbd} \frac{|u(x)|^p}{|x - y|^{k + sp}} \dd{s(x)} \dd{s(y)}.
    \end{align*}
    Note that
    \begin{align*}
        \int_{\manifold \setminus\genbd} \int_\genbd \frac{|u(x)|^p}{|x -y|^{k + sp}} \dd{s(x)} \dd{s(y)} = \int_{\manifold\setminus\genbd} \int_K \frac{|u(x)|^p}{|x -y|^{k + sp}} \dd{s(x)} \dd{s(y)}
        \\
        \leq \frac{|\manifold \setminus \genbd|}{\dist(\manifold \setminus \genbd, K)^{k + sp}} \Lpnorm{u}{\genbd}^p.
    \end{align*}
\end{proof}

\begin{lemma}
\label{lemma:extension at boundaries}
    Let $M$ be a $k$-dimensional Lipschitz submanifold of $\R^n$ and $\genbd \subset M$ open and connected. Furthermore, assume that $\partial \genbd$ is a $(k -1)$-dimensional Lipschitz submanifold of $M$. Finally, let $x \in \partial \genbd$ with $B \coloneqq \ball{\epsilon}{x}$ and $t_{x}: U_{x} \to B$ as defined in \cref{lem:local parameterizations}.

    Then, for every $u \in \Wsp{B \cap \genbd}$ with $s \in (0,1)$ and $p \in [1, +\infty)$, there exists an extension $\Bar{u}$ to $B \cap \manifold$ such that $\Bar{u}(y) = u(y)$ almost everywhere in $B \cap \genbd$ and $\Bar{u} \in \Wsp{B \cap \manifold}$ with
    \begin{align*}
        \Lpnorm{\Bar{u}}{B \cap \manifold}^p \leq C_{\partial}^{(0)} \Lpnorm{u}{B \cap \genbd}^p,
    && \text{ and } &&
        \Wspseminorm{\Bar{u}}{B \cap \manifold}^p \leq C_{\partial}^{(s)} \Wspseminorm{u}{B \cap \genbd}^p,
    \end{align*}
    where
    \begin{align*}
        C_{\partial}^{(0)} \coloneqq C_{x}^p IJ, &&
         C_{\partial}^{(s)} \coloneqq C_{x}^p (IJ)^2 (KL)^{k + sp},
    \end{align*}
    and $C_{x} > 0$ only depends on the pullback of $B \cap \partial \genbd$ to the chart $U_{x}$.
\end{lemma}

\begin{proof}
    First, we define $U_x^{-} \coloneqq t_x^{-1}(B \cap \genbd) \subset U_x$ and $v: U_x^- \to \R$ by setting $v(\hat{y}) \coloneqq u(t_x(\hat{y}))$ for $\hat{y} \in U_x^-$. We claim that $v \in \Wsp{U_x^-}$.
    Indeed, we have
    \begin{align*}
        \Lpnorm{v}{U_x^-}^p  = \int_{U_x^-} |u(t_x(\hat{x})|^p\dd{\hat{x} }
        = \int_{B \cap \genbd} |u(x)|^p J_{t_x^{-1}}(x) \dd{s(x)}
        \\
        \leq ||J_{t_x^{-1}}||_{L^\infty(B \cap \genbd
        )} \Lpnorm{u}{B \cap \genbd}^p \leq I \Lpnorm{u}{B \cap \genbd}^p,
    \end{align*}
    where $I$ is the bound from \cref{lem:local parameterizations}.
    Moreover, note that we have for $y, z \in B$ that
    \begin{align*}
        |y - z| = |t_x(t_x^{-1}(y)) - t_x(t_x^{-1}(z))| \leq L_{t_x} |t_x^{-1}(y) - t_x^{-1}(z)|.
    \end{align*}
    Therefore, we find for the seminorm,
    \begin{align*}
        \Wspseminorm{v}{U_x^-}^p &= \int_{U_x^-} \int_{U_x^-} \frac{|v(\hat{y}) - v(\hat{z})|^p}{|\hat{y} - \hat{z}|^{k + sp}} \dd{\hat{y}} \dd{\hat{z}} \\
        &= \int_{B \cap \genbd} \int_{B \cap \genbd} \frac{|u(y) - u(z)|^p}{|t_x^{-1}(y) - t_x^{-1}(z)|^{k + sp}} J_{t_x^{-1}}(y) J_{t_x^{-1}}(z) \dd{s(y)}\dd{s(z)}
        \\
        &\leq L_{t_x}^{k + sp} ||J_{t_x^{-1}}||^2_{L^\infty(B \cap \genbd)} \Wspseminorm{u}{B \cap \genbd}^p \leq L^{k + sp} I^2 \Wspseminorm{u}{B \cap \genbd}^p,
    \end{align*}
    where $L$ is the bound from \cref{lem:local parameterizations}. This proves the claim.
    
    Define $\eta \coloneqq B \cap \partial \genbd$. According to \cref{def:lipschitz domain}, $\partial \genbd$ is a $(k-1)$-dimensional Lipschitz boundary. Therefore, $U_x^-$ is a subset of $U_x$ with the Lipschitz boundary $\hat{\eta} \coloneqq t_x^{-1}(\eta)$, as $t_x$ is bi-Lipschitz by definition.
    By the Sobolev extension theorem for Lipschitz domains, there exists an extension $\Bar{v} \in \Wsp{U_x}$ of $v$ such that $\Bar{v}(\hat{y}) = v(\hat{y})$ almost everywhere in $U_x^-$ and
    \begin{align*}
        \Wspnorm{\Bar{v}}{U_x} \leq C_x \Wspnorm{v}{U_x^-},
    \end{align*}
    for some constant $C_x > 0$. We now define $\Bar{u}$ as the push-forward of $\Bar{v}$ by setting $\Bar{u}(x) \coloneqq \Bar{v}(t_x^{-1}(x)), x \in B \cap \manifold$. By the same argument as for $v$, we find
    \begin{align*}
        \Lpnorm{\Bar{u}}{B \cap \manifold}^p \leq J \Lpnorm{\Bar{v}}{U_x}^p \leq C_x^p IJ \Wspnorm{u}{B \cap \genbd}^p,
    \end{align*}
    and
    \begin{align*}
        \Wspseminorm{\Bar{u}}{B \cap \manifold}^p &\leq        
        J^2 K^{k + sp} \Wspseminorm{\Bar{v}}{U_x}^p \leq C_x^p (IJ)^2 (KL)^{k + sp}\Wspseminorm{u}{B \cap \genbd}^p,
    \end{align*}
    where $I, J, K$, and $L$ are the constants from \cref{lem:local parameterizations}.
\end{proof}

\begin{lemma}
\label{lemma:truncation on boundaries}
    Let $\manifold$ be a $k$-dimensional Lipschitz submanifold of $\R^n$ with a covering radius $\epsilon > 0$ as in \cref{lem:local parameterizations} and $\genbd \subset \manifold$ open and connected. Furthermore, let $u \in \Wsp{\genbd}$, with $s \in (0,1)$ and $p \in [1, +\infty)$, and $\psi$ a Lipschitz continuous function on $\genbd$ with $0 \leq \psi \leq 1$.

    Then, the product $\psi u$ is an element of $\Wsp{\genbd}$ with
    \begin{align*}
        \Lpnorm{\psi u}{\genbd} \leq \Lpnorm{u}{\genbd} 
        && \text{ and } &&
        \Wspseminorm{\psi u}{\genbd}^p \leq C_\psi^{(0)} \Lpnorm{u}{\genbd}^p +  C_\psi^{(s)}\Wspseminorm{u}{\genbd}^p,
    \end{align*}
    with $C_\psi^{(0)} \coloneqq 2^{p-1}$ and
    \begin{align*}
    C_\psi^{(s)} \coloneqq 2^{p-1} \left ( \frac{|\genbd|}{\epsilon^{k + (s-1)p}} + \frac{|\omega_k|}{p(1-s)} \frac{JK^k}{\epsilon^{(s-1)p}} \right),
    \end{align*}
    where $|\omega_k|$ denotes the surface area of the $k$-dimensional unit sphere.
\end{lemma}

\begin{proof}
    We clearly have $\Lpnorm{\psi u}{\genbd} \leq \Lpnorm{u}{\genbd}$ by monotonicity. By adding and subtracting the factor $\psi(x)u(y)$ in the numerator, we find
    \begin{align*}
        \Wspseminorm{\psi u}{\genbd}^p \leq 2^{p-1}\left ( 
            \int_\genbd \int_\genbd \frac{|\psi(x)u(x) - \psi(x) u(y)|^p}{|x - y|^{k + sp}} \dd{s(x)} \dd{s(y)} \right.\nonumber\\
            \left.+ \int_\genbd \int_\genbd \frac{|\psi(x)u(y) - \psi(y) u(y)|^p}{|x - y|^{k + sp}} \dd{s(x)} \dd{s(y)}
        \right)\\
        \leq 2^{p-1} \left( \Wspseminorm{u}{\genbd}^p + \int_\genbd \int_\genbd \frac{|u(x)|^p |\psi(x) - \psi(y)|^p}{|x - y|^{k + sp}} \dd{s(x)} \dd{s(y)}\right).
    \end{align*}
    Let $\epsilon > 0$ be as in \cref{lem:local parameterizations} and $L_\psi$ the Lipschitz constant of $\psi$. Then, interchanging of the integrals according to Fubini's Theorem in the second integral yields,
    \begin{align}
        \label{eq:truncation split estimate}
        \int_\genbd \int_\genbd \frac{|u(x)|^p |\psi(x) - \psi(y)|^p}{|x - y|^{k + sp}} \dd{s(y)} \dd{s(x)}
        \nonumber \\
        \leq L_\psi^p \left(
        \int_\genbd \int_{\genbd \cap |x - y| < \epsilon} \frac{|u(x)|^p}{|x - y|^{k + (s-1)p}} \dd{s(y)} \dd{s(x)} \right.
        \nonumber \\
        \left. + \int_\genbd \int_{\genbd \cap |x - y| \geq \epsilon} \frac{|u(x)|^p}{|x - y|^{k + (s-1)p}} \dd{s(y)} \dd{s(x)} 
        \right).
    \end{align}
    
    For the second term, we immediately find
    \begin{align*}
        \int_\genbd \int_{\genbd \cap |x - y| \geq \epsilon} \frac{|u(x)|^p}{|x - y|^{k + (s-1)p}} \dd{s(y)} \dd{s(x)}
        \leq \frac{|\genbd|}{\epsilon^{k + (s-1)p}} \Lpnorm{u}{\genbd}^p.
    \end{align*}

    Using the chart map $t_x: U_x \to \ball{\epsilon}{x} \cap \manifold$, we find together with the change of variable formula in \cref{eq:change of variable formula} that
    \begin{align*}
        \int_{\genbd \cap |x - y| < \epsilon} \frac{1}{|x - y|^{k + (s-1)p}} \dd{s(y)} = \int_{U_x} \frac{1}{|t_{x}(0) - t_{x}(\hat y)|^{k + (s-1)p}}J_{t_{x}}(\hat{y}) \dd{\hat{y}}
        \\
        \leq ||J_{t_{x}}||_{L^\infty(U_x)} L_{t_{x}^{-1}}^{k + (s-1)p} \int_{U_x} \frac{1}{|\hat{y}|^{k + (s-1)p}} \dd{\hat{y}}
        \\
        \leq J K^{k + (s-1)p} \int_{\ball{\epsilon K}{0}}  \frac{1}{|\hat{y}|^{k + (s-1)p}} \dd{\hat{y}},
    \end{align*}
    where we used assertion 3 in \cref{lem:local parameterizations}.
    
    Note that the kernel $r \mapsto r^{-(1 + (s-1)p)}$ is integrable on the interval $[0, \epsilon K]$ whenever $p(1 - s) > 0$. The latter is always the case for $p \geq 1$ and $s \in (0,1)$. Therefore, we have
    \begin{align*}
        \int_{\ball{\epsilon K}{0}}  \frac{1}{|\hat{y}|^{k + (s-1)p}} \dd{\hat{y}} = |\omega_k| \int_0^{\epsilon K} \frac{1}{r^{1 + (s-1)p}} \dd{r} = |\omega_k| \frac{(\epsilon K)^{p(1-s)}}{p(1-s)}.
    \end{align*}
    Reinserting this expression into the corresponding double integral in \cref{eq:truncation split estimate} and using the scaling properties in assertion 4 of \cref{lem:local parameterizations}, we find
    \begin{align*}
        \int_\genbd \int_{\genbd \cap |x - y| < \epsilon} \frac{|u(x)|^p}{|x - y|^{k + (s - 1)p}} \dd{s(y)} \dd{s(x)}
        \leq \frac{|\omega_k|}{p(1-s)} \epsilon^{p(1 - s)} JK^k\Lpnorm{u}{\genbd}^p.
    \end{align*}
    The claim then follows from gathering the above estimates.
\end{proof}

The constants in \cref{lemma:zero extension on boundaries,lemma:extension at boundaries,lemma:truncation on boundaries} depend on the geometry of $\manifold$ and $\partial \genbd$. Since these vary under rescaling, we derive their explicit dependence on $|\genbd|$ in \cref{lemma:scaling properties}.

\begin{lemma}
\label{lemma:scaling properties}
Let $\manifold$ be a $k$-dimensional Lipschitz submanifold of $\R^n$ and $\genbd \subset \manifold$ open and connected. Furthermore, let $\hat{\manifold}$ denote the manifold obtained after rescaling $\R^n$, such that $|\hat{\genbd}| = 1$. We recall the constants $C^\star, C_{\partial}^{(0)}, C_{\partial}^{(s)}, C_{\psi}^{(0)}$, and $C_{\psi}^{(s)}$ from \cref{lemma:zero extension on boundaries,lemma:extension at boundaries,lemma:truncation on boundaries} that are defined in terms of the geometric constants $I, J, K, L$ and $\epsilon$ from \cref{lem:local parameterizations}.
Then, the following scaling relations hold with respect to their rescaled counterparts denoted by $\hat{\cdot}$:
\begin{align*}
    C^\star &= \hat{C}^\star |\genbd|^{-\frac{sp}{k}}, &
    C^{(0)}_\partial &= \hat{C}^{(0)}_\partial, &
    C^{(s)}_\partial &= \hat{C}^{(s)}_\partial, &
    \\
    C_\psi^{(0)} &= \hat{C}_\psi^{(0)}|\genbd|^{\frac{(1 - s)p}{k}}, &
    C_\psi^{(s)} &= \hat{C}_\psi^{(s)}. &&
\end{align*}
\end{lemma}

\begin{proof}
To prove these scaling relations, we first establish the behavior of the uniform Lipschitz constants $L$ and $K$ and the uniform bounds on the Jacobian determinants $J$ and $I$.
To that end, we rescale $\manifold$ by a factor $\lambda > 0$ denoted by $\manifold_\lambda$. Then, for $x \in \manifold$ we have that $x_\lambda \coloneqq \lambda x \in \manifold_\lambda$. In particular, from a chart map $t_x: U_x \to \ball{\epsilon}{x} \cap \manifold$ for $\manifold$, we obtain a valid chart map $\Tilde{t}_{x_\lambda}: U_x \to \ball{\lambda \epsilon}{x_\lambda} \cap \manifold_\lambda$ for $x_\lambda \in \manifold_\lambda$ by setting $\Tilde{t}_{x_\lambda}(\hat{y}) \coloneqq \lambda t_x(\hat{y})$ for $\hat{y} \in U_x$. It follows immediately that
\begin{align*}
    L_{\Tilde{t}_{x_\lambda}} = \lambda L_{t_x} & &\text{ and } & & ||J_{\Tilde{t}_{x_\lambda}}||_{L^\infty(U_x)} = \lambda^k ||J_{t_{x}}||_{L^\infty(U_x)},
\end{align*}
and inversely for $L_{\Tilde{t}_{x_\lambda}^{-1}}$ and $||J_{\Tilde{t}_{x_\lambda}^{-1}}||_{L^\infty(\ball{\lambda\epsilon}{x} \cap \manifold_\lambda)}$. Rescaling to reference size $|\genbd_\lambda| = 1$ requires $\lambda = |\genbd|^{-\frac{1}{k}}$. Inserting this choice of $\lambda$ in the above expression and dividing by $\lambda$ and $\lambda^k$, respectively, yields
\begin{align}
\label{eq:uniform lipschitz constant:scaling}
    L &= |\genbd|^{-\frac{1}{k}} \hat{L}, & K &= |\genbd|^{\frac{1}{k}} \hat{K}, 
    \\
\label{eq:uniform jacobian constant:scaling}
    J &= |\genbd|^{-1} \hat{J}, & I &= |\genbd| \hat{I}.
\end{align}
Furthermore, the covering radius $\epsilon > 0$ from \cref{lem:local parameterizations} satisfies $\epsilon = \hat{\epsilon} |\genbd|^{\frac{1}{k}}$.

For $C^\star$, we observe that
\begin{align*}
    C^\star = \frac{|\manifold \setminus \genbd|}{\dist(\manifold \setminus \genbd, F)^{k + sp}} = \frac{|\hat{\manifold} \setminus \hat{\genbd}| |\genbd|}{\dist(\hat{\manifold} \setminus \hat{\genbd}, \hat{F})^{k + sp} |\genbd|^{1 + \frac{sp}{k}}} = \hat{C}^\star |\genbd|^{-\frac{sp}{k}},
\end{align*}
as claimed.

For $C^{(0)}_\partial$ and $C^{(s)}_\partial$, we note that the constant $C_x$ in the proof of \cref{lemma:extension at boundaries}, that stems from the extension in the local chart $U_x$, is invariant under rescaling as the latter is accounted for by the chart map $t_x$. Hence, we conclude that $C_x = \hat{C}_x$. Using the scaling relations from \cref{eq:uniform lipschitz constant:scaling,eq:uniform jacobian constant:scaling}, we thus find
\begin{align*}
    C^{(0)}_\partial = C_x^p IJ = \hat{C}_x^p \hat{I}\hat{J} = \hat{C}^{(0)}_\partial,
\end{align*}
and similarly for $C^{(s)}_\partial$.

Clearly, $C^{(s)}_\psi = \hat{C}^{(s)}_\psi$. Furthermore, we have
\begin{align*}
    \frac{|\genbd|}{\epsilon^{k + (s - 1)p}} = \frac{|\hat{\genbd}|}{\hat{\epsilon}^{k + (s - 1)p}} |\genbd|^{\frac{(1 - s)p}{k}}.
\end{align*}
Using again the scaling relations from \cref{eq:uniform lipschitz constant:scaling,eq:uniform jacobian constant:scaling}, we obtain
\begin{align*}
    \frac{J K^k}{\epsilon^{(s-1)p}} = \frac{\hat{J} \hat{K}}{\hat{\epsilon}^{(s - 1)p}} |\genbd|^{-\frac{(1 - s)p}{k}}.
\end{align*}
Therefore, we conclude that $C^{(0)}_\psi = \hat{C}^{(0)}_\psi |\genbd|^{\frac{(1 - s)p}{k}}$, as claimed.
\end{proof}

With the above results, we are ready to prove the main theorem.

\begin{proof}[Proof of \cref{thm:boundary extension operator}]
    Let $\epsilon > 0$ be as in \cref{lem:local parameterizations}. Since $\partial \genbd$ is compact, we find finitely many points $x_1, \ldots, x_N \in \partial \genbd$ such that $\partial \genbd \subset \bigcup_{i = 1}^N B_i$, where $B_i \coloneqq \ball{\epsilon}{x_i} \cap \manifold, i = 1, \ldots, N$. Let $B_0 \subset \genbd$ be another open subset such that the collection $\{B_i\}_{i=0}^N$ is a finite open covering of the closure $\overline{\genbd}$. Let $\{\psi_i\}_{i=0}^N$ be a partition of unity subordinated to the latter covering. By definition of the partition of unity, we have
    \begin{align*}
        u = \sum_{i = 0}^N \psi_i \, u \; \text{ in } \genbd.
    \end{align*}

    Let us first consider the term $\psi_0 u$. By \cref{lemma:truncation on boundaries}, $\psi_0 u$ is an element of $\Wsp{\genbd}$ and vanishes in a neighborhood of $\partial \genbd$. Hence, by \cref{lemma:zero extension on boundaries} and \cref{lemma:truncation on boundaries}, the extension by zero $w_0 \coloneqq (\psi_0 u)^\star$ to $\manifold$ is an element of $\Wsp{\manifold}$ with
    \begin{nalign}
    \label{eq:proof:zero term estimate}
        \Wspnorm{w_0}{\manifold}^p \leq \left(1 + C^\star_0\right) \Lpnorm{\psi_0 u}{\genbd}^p + \Wspseminorm{\psi_0 u}{\genbd}^p
        \\
        \leq \left(1 + C^\star_0 + C_{\psi_0}^{(0)}\right) \Lpnorm{u}{\genbd}^p + C_{\psi_0}^{(s)}\Wspseminorm{u}{\genbd}^p, 
    \end{nalign}
    where $C^\star_0, C_{\psi_0}^{(0)}$, and $C_{\psi_0}^{(s)}$ denote the constants from \cref{lemma:zero extension on boundaries,lemma:truncation on boundaries} corresponding to $B_0$.
    After applying the scaling relations from \cref{lemma:scaling properties}, \cref{eq:proof:zero term estimate}, we find constants $\Tilde{C}_0^{(0)} > 0$ and $\Tilde{C}_0^{(s)} > 0$ independent of the measure $|\genbd|$, such that
    \begin{align}
        \Wspnorm{w_0}{\manifold}^p
        \leq \Tilde{C}_0^{(0)} \left(1 + |\genbd|^{-\frac{sp}{k}} + |\genbd|^{\frac{(1 - s)p}{k}}\right) \Lpnorm{u}{\genbd}^p + \Tilde{C}_0^{(s)}\Wspseminorm{u}{\genbd}^p. 
    \end{align}

    Next, we treat the terms intersecting the boundary $\partial \genbd$. For $1 \leq i \leq N$, we have $u_{|B_i \cap \genbd} \in \Wsp{B_i\cap \genbd}$. By \cref{lemma:extension at boundaries}, there exists an extension $v_i \in \Wsp{B_i \cap \manifold}$ with
    \begin{align}
    \label{eq:symmetric extension lp}
        \Lpnorm{v_i}{B_i \cap \manifold}^p \leq C_{\partial_i}^{(0)} \Lpnorm{u_{|B_i \cap \genbd}}{B_i \cap \genbd}^p \leq C_{\partial_i}^{(0)} \Lpnorm{u}{\genbd}^p,
    \end{align}
    and
    \begin{align}
    \label{eq:symmetric extension seminorm}
        \Wspseminorm{v_i}{B_i \cap \manifold}^p \leq C_{\partial_i}^{(s)} \Wspseminorm{u_{|B_i \cap \genbd}}{B_i \cap \genbd}^p \leq C_{\partial_i}^{(s)} \Wspseminorm{u}{\genbd}^p,
    \end{align}
    where we write $C_{\partial_i}^{(0)}$ and $C_{\partial_i}^{(s)}$ for the constants in \cref{lemma:extension at boundaries} corresponding to the ball $B_i$.
    In particular, we must have $\psi_i v_i = \psi_i u$ almost everywhere in $B_i \cap \genbd$. By definition, $\psi_i v_i$ has compact support in $B_i \cap \manifold$. Therefore, we can apply once more  \cref{lemma:zero extension on boundaries} and define the extension by zero $w_i \coloneqq (\psi_iv_i)^\star \in \Wsp{\manifold}$, such that by \cref{lemma:zero extension on boundaries}, \cref{lemma:truncation on boundaries}, \cref{eq:symmetric extension lp,eq:symmetric extension seminorm}, we have
    \begin{nalign}
    \label{eq:proof:boundary terms estimate}
        \Wspnorm{w_i}{\manifold}^p \leq \left(1 + C^\star_i + C_{\psi_i}^{(0)}\right) \Lpnorm{v_i}{B_i \cap \manifold}^p + C_{\psi_i}^{(s)}\Wspseminorm{v_i}{B_i \cap \manifold}^p
        \\
        \leq \left( 1 + C^\star_i + C_{\psi_i}^{(0)} \right) C_{\partial_i}^{(0)} \Lpnorm{u}{\genbd}^p +  C_{\partial_i}^{(s)} C_{\psi_i}^{(s)}\Wspseminorm{u}{\genbd}^p,
    \end{nalign}
    where $C^\star_i, C_{\psi_i}^{(0)}$, and $C_{\psi_i}^{(s)}$ denote the constants from \cref{lemma:zero extension on boundaries,lemma:truncation on boundaries} corresponding to the ball $B_i$.
    Applying once more the scaling relationships from \cref{lemma:scaling properties}, \cref{eq:proof:boundary terms estimate} reads
    \begin{align*}
        \Wspnorm{w_i}{\manifold}^p \leq \Tilde{C}_i^{(0)} \left(1 + |\genbd|^{-\frac{sp}{k}} + |\genbd|^{\frac{(1 - s)p}{k}}\right) \Lpnorm{u}{\genbd}^p + \Tilde{C}_i^{(s)}\Wspseminorm{u}{\genbd}^p,
    \end{align*}
    for constants $\Tilde{C}_i^{(0)} > 0$ and $\Tilde{C}_i^{(s)}$ independent of the measure $|\genbd|$.

    Defining the extension $\extop{\genbd}{\manifold}(u)$ by
    \begin{align*}
        \extop{\genbd}{\manifold}(u)(x) \coloneqq \sum_{i = 0}^N w_i(x), \quad x \in \manifold,
    \end{align*}
    we have by construction that $\extop{\genbd}{\manifold}(u) = u$ almost everywhere in $\genbd$ and we find constants $C^{(0)} > 0$ and $C_i^{(s)} > 0$ independent of the measure $|\genbd|$, such that
    \begin{align*}
        \Wspnorm{\extop{\genbd}{\manifold}(u)}{\manifold}^p &\leq \sum_{i = 0}^N \Wspnorm{w_i}{\manifold}^p
        \\
        &\leq
        C^{(0)} \left(1 + |\genbd|^{-\frac{sp}{k}} + |\genbd|^{\frac{(1-s)p}{k}}\right) \Lpnorm{u}{\genbd}^p + C^{(s)} \Wspseminorm{u}{\genbd}^p.
    \end{align*}
    Since the decomposition by partition of unity, the local pullbacks and extensions are all linear operations, the resulting extension operator is linear too.
    This concludes the proof.
\end{proof}

\begin{remark}
\label{remark:generalization}
The extension result of \cref{thm:boundary extension operator} can in principle be generalized to Sobolev exponents $s \geq 1$. In this case, however, one must assume higher regularity of the manifold and of $\partial\Lambda$, since composition with local charts preserves $W^{s,p}$ only for sufficiently smooth parametrizations.

Sobolev spaces of order $s \geq 1$ on Lipschitz or smooth manifolds can be defined either via local coordinate charts, see \cite{grisvard_elliptic_2011}, or intrinsically by means of a Riemannian metric and covariant derivatives, see \cite{hebey_nonlinear_2000}. While these constructions yield equivalent norms, the associated constants depend on higher-order geometric quantities such as derivatives of the charts or curvature bounds.

As a consequence, although extension operators continue to exist in this setting, obtaining explicit control of their continuity constants under rescaling becomes substantially more involved. Therefore, we restrict our analysis to $s \in (0,1)$.
\end{remark}

\section{Conclusion}

In this paper, we have established an extension theorem for fractional Sobolev spaces $W^{s,p}$ on subsets $\genbd$ of bounded Lipschitz submanifolds, provided that $\genbd$ itself has a Lipschitz boundary relative to the manifold. In particular, the dependence on the size of the subset is explicit in all estimates.
Our proof adapts the classical arguments for Lipschitz domains in Euclidean space from \cite{di_nezza_hitchhikers_2012}, with particular attention paid to the uniform control of local parameterizations on the compact Lipschitz manifold $\manifold$, as detailed in \cref{lem:local parameterizations}.

Future work includes generalizing the result to Sobolev exponents $s \geq 1$ as outlined in \cref{remark:generalization} under stronger regularity assumptions on the manifold.

\subsection*{Acknowledgments}
The author gratefully acknowledges Annalisa Buffa for her insightful comments and discussions.

\bibliography{references}

\subsection*{Statements and declarations}

\subsubsection*{Funding}
This research was supported by the Swiss National Science Foundation via project MINT n. 200021\_215099, PDE tools for analysis-aware geometry processing in simulation science.

\subsubsection*{Competing interests}
The author has no relevant financial or non-financial interests to disclose.

\subsubsection*{Data availability}
This manuscript has no associated data.

\end{document}